\documentclass[11pt,letterpaper]{article}
\usepackage{amsmath,amssymb,amsfonts,mathrsfs}
\usepackage{graphicx}
\usepackage[margin=1in]{geometry}

\usepackage{amssymb}
\usepackage{graphicx,longtable,wrapfig}
\usepackage{amsmath}
\usepackage{epstopdf}

\usepackage{url}

\usepackage{caption,subcaption}

\newtheorem{theorem}{Theorem}
 
\title{Assortativity in generalized preferential attachment models\thanks{This is a journal version of the paper appeared in Proc. WAW'16, LNCS 10088, pp. 9-21, 2016.}
	\thanks{This work is supported by the Russian President grant MK-527.2017.1.}}
 
\author{Alexander Krot$^1$ \and Liudmila Ostroumova Prokhorenkova$^{1,2}$}
\date{%
    $^1$Moscow Institute of Physics and Technology, Moscow, Russia\\%
    $^2$Yandex, Moscow, Russia\\[2ex]%
}

\begin{document}
\maketitle
\newcommand{\E}{\mathbb{E}}
\newcommand{\Prob}{\mathsf{P}}
\newcommand{\G}{\mathsf{\Gamma}}
\newcommand{\BA}{Barab\'asi--Albert}
\newcommand{\BO}{Buckley--Osthus}
\newcommand{\HK}{Holme and Kim}

\newcommand{\widesim}[2][1.5]{
  \mathrel{\overset{#2}{\scalebox{#1}[1]{$\sim$}}}
}

\begin{abstract}
In this paper, we analyze assortativity in a wide class of preferential attachment models (PA-class). It was previously shown that the degree distribution in all models of the PA-class follows a power law. Also, the global and the average local clustering coefficients were analyzed. We expand these results by analyzing the assortativity property of the PA-class of models. Namely, we analyze the behavior of $d_{nn}(d)$ which is the average degree of a neighbor of a vertex of degree~$d$.


\end{abstract}

\section{Introduction}\label{Introduction}

Nowadays, there is a great deal of interest in structure and dynamics of real-world networks, from Internet and society networks~\cite{BA_Review,Networks,Chayes} to biological networks~\cite{BioInfoPrior}. The key problem is how to build a model which describes the properties of a given network. Such models are used in physics, information retrieval, data mining, bioinformatics, etc.~\cite{BA_Review,Networks,Math_Results,Lescovec}.

Real-world networks have some common properties~\cite{Girvan,Newman2,Newman3,Watts}.
For example, for the majority of studied networks, the degree distribution was observed to follow the power law, which means that the portion of vertices with degree $d$ decreases as $d^{-\gamma}$ for some $\gamma>0$~\cite{BA,Networks,Broder,F-F-F}. Another important property of complex networks is high clustering coefficient~\cite{Newman3} which, roughly speaking, measures how likely two neighbors of a vertex are connected.

Another key metric in complex networks analysis is the assortativity coefficient which was first introduced by Newman~\cite{Newman} as the Pearson's correlation coefficient for the pairs $\{(d_i, d_j) | e_{ij} \in E \}$. In assortative graphs edges tend to connect vertices of similar degrees, while in disassortative networks vertices of low degree are more likely to be connected to vertices of high degree. Assortativity coefficient lies between -1 and 1; when this coefficient equals 1, the network is said to have perfect assortative mixing patterns, when it equals 0, the network is non-assortative, while at -1 the network is completely disassortative.
However, as discussed in \cite{DDC,DDD}, despite Pearson's correlation coefficient is most commonly used to measure assortativity of networks, this coefficient is size-depend when the degree distribution has infinite variance.

Another way to analyze assortativity is to consider the behavior of $d_{nn}(d)$~--- the average degree of a neighbor of a vertex of degree $d$~--- introduced in~\cite{dnn}. A graph is called assortative if $d_{nn}(d)$ is an increasing function of $d$, whereas it is referred to as disassortative when $d_{nn}(d)$ is a decreasing function of $d$.
We analyze $d_{nn}(d)$ instead of measuring the correlation since the obtained function of $d$ can give a deeper insight into the network structure.

It was previously shown that in some real-world networks $d_{nn}(d)$ behaves as $d^\nu$ for some $\nu$, which can be positive (assortative networks) or negative (disassortative networks)~\cite{Networks,Dist_D}. Interestingly, as we show in this paper, in a wide class of preferential attachment models $\E d_{nn}(d) \propto \log(d)$ as $d \to \infty$.

Assortativity has many applications, for instance, it can be used in the epidemiology. In social networks we usually observe assortative mixing, so diseases targeting high degree individuals are likely to spread to other high degree nodes. On the other hand, biological networks are usually disassortative, therefore vaccination strategies that specifically target the high degree vertices may quickly destroy the epidemic network.

In this paper, we study the behavior of $d_{nn}(d)$ in the T-subclass of the PA-class of models, which was first introduced in~\cite{GPA}. This class includes many well-known models based on the preferential attachment principle: LCD~\cite{LCD_degrees}, Buckley-Osthus~\cite{Buckley_Osthus}, Holme-Kim~\cite{Holme_Kim}, RAN~\cite{RAN}, etc. Despite the fact that the T-subclass generalizes many different models, we are able to rigorously analyze $d_{nn}(d)$ in the whole class of models for $\gamma>3$ (the case of a finite variance). In particular, we prove that in this case the expectation of $d_{nn}(d)$ asymptotically behaves as $\log(d)$ (up to a constant multiplier).
However, this approximation works reasonably well only for very large values of $d$ and for $d<10^4$ we observe a different behavior which may look like $d^\nu$ for some $\nu>0$, as it was observed in some real-world networks.

However, if the degree distribution has infinite variance, $d_{nn}(d)$ is much harder to analyze. For the configuration model this problem is addressed in~\cite{ANND}. Namely, it is shown that when the variance of the degree distribution is infinite, $d_{nn}(d)$ scales with $n$, where $n$ is the number of vertices, and a corresponding central limit theorem is proven. 
It turns out that similar problem arises in the PA-class of models, and in the case of infinite variance the precise asymptotics for $d_{nn}(d)$ cannot be obtained, as we discuss further in the paper. However, using some heuristics, we are able to make hypotheses on the expected behavior of $d_{nn}(d)$. Thus, for $\gamma = 3$ we obtain $d_{nn}(d) \propto \log(n)$ and for $\gamma<3$ the asymptotics is $d_{nn}(d) \propto d^{-(3-\gamma)}n^{\frac{3-\gamma}{\gamma-1}}$. It turns out that these hypotheses agree well with our simulations.

The remainder of the paper is organized as follows.
In Section~\ref{PAclass}, we give a formal definition of the PA-class and
present some relevant known results.
In Section~\ref{Assort}, we state new theoretical results on the behavior of $d_{nn}(d)$; then we prove all the theorems; and finally we discuss our hypotheses for the case $\gamma \le 3$. In Section~\ref{experiments}, we make some simulations in order to illustrate our results for $d_{nn}(d)$. Section~\ref{sec:conclusion} concludes the paper.


\section{Generalized Preferential Attachment}\label{PAclass}

\subsection{Definition of the PA-class}\label{class}

Let us formally define the PA-class of models which was first proposed in \cite{GPA}.
Let $G_{m}^n$ ($m$ is a constant parameter of the model, $n \ge n_0$) be a graph with $n$ vertices $\{1, \ldots, n\}$ and $mn$ edges obtained as a result of the following process. We start at the time $n_0$ from an arbitrary graph $G_{m}^{n_0}$ with $n_0$ vertices and $m n_0$ edges. On the $(n+1)$-th step ($n\geq n_0$), we make the graph $G_{m}^{n+1}$ from $G_{m}^{n}$ by adding a new vertex $n+1$ and $m$ edges connecting this vertex to some $m$ vertices from the set $\{1, \ldots , n, n+1\}$. Denote by $d_v^{n}$ the degree of a vertex $v$ in~$G_m^n$. If for some constants $A$ and $B$ the following conditions are satisfied
\begin{equation}\label{OneStepChangedDegreeDraft}
\Prob\left( d_v^{n+1} = d_v^{n} \mid G_m^{n}\right) = 1 - A \frac{d_v^{n}}{n} - B\frac{1}{n} + O\left(\frac{\left(d_v^n\right)^2}{n^2}\right),\,\,1 \le v \le n \;,
\end{equation}
\begin{equation}\label{OneStepChangedDegreeDraft_2}
\Prob\left( d_v^{n+1} = d_v^{n} + 1 \mid G_m^{n}\right) =  A \frac{d_v^{n}}{n} + B\frac{1}{n} + O\left(\frac{\left(d_v^n\right)^2}{n^2}\right), \,\,1 \le v \le n \;,
\end{equation}
\begin{equation}\label{OneStepChangedDegreeDraft_3}
\Prob\left( d_v^{n+1} = d_v^{n} + j \mid G_m^{n}\right) =  O\left(\frac{\left(d_v^n\right)^2}{n^2}\right), \,\,
2\le j \le m,\,\,1 \le v \le n \;,
\end{equation}
\begin{equation}\label{LoopProbability}
\Prob( d_{n+1}^{n+1} =  m + j ) = O\left(\frac 1 n \right), \,\,
1\le j \le m\;,
\end{equation}
then the random graph process $G_m^n$ is a model from the PA-class. Note that $O(\cdot)$ above are defined as $n\to\infty$.
Here, as in~\cite{GPA}, we require $2mA + B = m$ and $0 \le A \le 1$.
We further omit $n$ in $d_j^n$ for simplicity of notation.

As it is explained in \cite{GPA}, even fixing values of parameters $A$ and $m$ does not specify a concrete procedure for constructing a network.
There are a lot of models possessing very different properties and satisfying the conditions~(\ref{OneStepChangedDegreeDraft}--\ref{LoopProbability}), e.g., LCD~\cite{LCD_degrees}, Buckley-Osthus~\cite{Buckley_Osthus}, Holme-Kim~\cite{Holme_Kim}, and RAN~\cite{RAN} models.

\subsection{Power-law Degree Distribution}\label{DegreeDistribution}

Let $N_n(d)$ be the number of vertices of degree $d$ in~$G_m^n$.
The following theorems on the expectation of $N_n(d)$ and its concentration were proved in~\cite{GPA}.

\begin{theorem}\label{Expectation}
For every model in PA-class and for every $d=d(n) \ge m$
$$
\E N_n(d) = c(m,d) \left(n + O\left(d^{2 + \frac{1}{A}}\right)\right)\,,
$$
where
\begin{equation}\label{Constant}
c(m,d) = \frac{\G\left(d + \frac{B}{A}\right)\G\left(m + \frac{B+1}{A}\right)}{A\,\G\left(d + \frac{B+A+1}{A}\right)\G\left(m + \frac{B}{A}\right)}
\widesim{d \rightarrow \infty} \frac{\G\left(m + \frac{B+1}{A}\right)d^{-1-\frac{1}{A}}}{A\, \G\left(m + \frac{B}{A}\right)}
\end{equation}
and $\G(x)$ is the gamma function.
\end{theorem}

\begin{theorem}\label{Concentration}
For every model from the PA-class and for every $d=d(n)$ we have
$$
\Prob\left(|N_n(d) - \E N_n(d)| \ge d \, \sqrt{n}  \, \log{n}\right) = O\left(n^{-\log{n}}\right).
$$
\end{theorem}
These two theorems mean that the degree distribution follows (asymptotically) the power law with the parameter $\gamma = 1+\frac{1}{A}$.

\subsection{Clustering Coefficient}\label{ClusteringCoefficient}

A T-subclass of the PA-class was introduced in~\cite{GPA}.
In this case, the following additional condition is required:
\begin{equation}\label{D_definition}
\Prob\left( d_i^{n+1} = d_i^{n} + 1, d_j^{n+1} = d_j^{n} + 1 \mid G_m^{n}\right) = e_{ij} \frac{D}{mn} + O\left(\frac{d_i^{n} d_j^{n}}{n^2}\right) \;,
\end{equation}
where $1 \le i,j \le n$, $e_{ij}$ is the number of edges between the vertices $i$ and $j$ in $G_m^n$ and $D$ is a non-negative constant. Note that this property still does not define the correlation between edges completely, but it is sufficient for studying the clustering coefficients. Also, this subclass still covers all well-known models mentioned above.

There are two well-known definitions of the clustering coefficient of a graph $G$.
The \emph{global clustering coefficient} $C_1(G)$ is the ratio of three times the number of triangles to the number of pairs of adjacent edges in $G$. The \emph{average local clustering coefficient} is defined as $C_2(G) = \frac{1}{n} \sum_{i=1}^n C^i$, where $C^i$ is the local clustering coefficient for a vertex $i$: $C^i= \frac{T^i}{P_2^i}$, $T^i$ is the number of edges between the neighbors of the vertex $i$ and $P_2^i$ is the number of pairs of neighbors.

The clustering coefficients for the T-subclass were analyzed in~\cite{GPA} and~\cite{LCC}. For example, in~\cite{GPA} it was proven that in some cases ($2A\ge1$) the global clustering coefficient $C_1(G_m^n)$ tends to zero as the number of vertices grows for all models from the PA-class. Additionally, it was shown that the average local clustering coefficient $C_2(G_m^n)$ does not tend to zero for the T-subclass with $D>0$.
In \cite{LCC} the local clustering coefficient averaged over the vertices of degree $d$ was analyzed. It was proven that this coefficient $C(d)$ asymptotically decreases as $\frac{2D}{Am} \cdot d^{-1}$ for $A<\frac{3}{4}$.

\section{Theoretical results}\label{Assort}

In this section, we first present our theoretical results on the behavior of $d_{nn}(d)$ in the T-subclass with $A<1/2$. Second, we prove these results. Finally, we discuss our hypotheses about the behavior of $d_{nn}(d)$ for the case $A\ge \frac{1}{2}$.

\subsection{Main results ($A < 1/2$)}\label{main_results}

Denote by $S_n(d)$ the sum of the degrees of all neighbors of all vertices of degree $d$: 
$$
S_n(d)=\sum_{i:d_i=d}\,{\sum_{j:ij \in E(G)}{d_j}}\,,
$$
where $E(G)$ is the set of edges of the graph $G$.
One possible way to analyze the assortativity of an undirected graph $G$ is to consider the average degree of the neighbors of vertices with a given degree $d$:
\begin{equation}\label{CD_definition}
d_{nn}(d) = \frac{S_n(d)}{N_n(d) \cdot d}\,.
\end{equation}
If $d_{nn}(d)$ is an increasing function of $d$, then the network is called assortative.
Vice-versa, in the disassortative case $d_{nn}(d)$ decreases.

In order to estimate $\E d_{nn}(d)$, we first estimate $\E S_n(d)$ and then use Theorems~\ref{Expectation} and~\ref{Concentration} on the behavior of $N_n(d)$.
Namely, we prove the following theorems.

\begin{theorem}\label{PASexp}
Let $G_m^n$ belong to the T-subclass with $A<\frac 1 2$. Then, for any $\varepsilon>0$ and every $d=d(n)\ge m$
\begin{equation}\label{eq:main}
\E S_n(d) = M(d)\left(n + O \left(n^{2A+\varepsilon} d^{\xi}\right)\right),
\end{equation}
where $\xi = \max\left\{3+\frac 1 A-4A,\frac{2}{1-A}\right\}$ and
\begin{equation}\label{eq:M}
M(d) = \left(Ad+B+1\right) \left(\frac{X}{Am+B+1} + \sum_{i=m+1}^{d} {Y(i)}\right) \cdot c(m,d) ,
\end{equation}
\begin{equation*}
X = \frac{m}{A(m-1)+B+1} \left(B-\frac{D}{m}+\frac{\left(A(m-1)+2B+1\right) \cdot (Am+B+1)}{1-2A} \right),
\end{equation*}
\begin{equation*}
Y(i) = \frac{1}{A(i-1)+B+1} \left(\frac{(B-D/m)i}{Ai+B+1}+\frac{(D/m)\cdot(i-1)}{A(i-1)+B} +m \right).
\end{equation*}
Asymptotically we have
\begin{equation}\label{eq:M_limit}
M(d)\widesim{d \rightarrow \infty} \frac{Am+B}{A^2} \cdot \frac{\G\left(m + \frac{B+1}{A}\right)}{\G\left(m + \frac{B}{A}\right)} \cdot \log(d) \cdot d^{-\frac{1}{A}}.
\end{equation}
\end{theorem}

\begin{theorem}\label{PACD}
Let $G_m^n$ belong to the T-subclass of the PA-class with $A<\frac 1 2$. Then for any $\varepsilon>0$ and for every $d=d(n)\ge m$
$$
\E d_{nn}(d) = \frac{M(d)}{d \, c(m,d)}\left(1 + O\left(\frac{n^{2A+\varepsilon} d^{\xi}}{n} + \frac{d^{2+\frac 1 A} \log{n}}{\sqrt{n}} \right) \right),
$$
where $\xi = \max\left\{3+\frac 1 A-4A,\frac{2}{1-A}\right\}$.
Note that asymptotically 
\begin{equation}\label{eq:dnn_limit}
\frac{M(d)}{d \cdot c(m,d)} \widesim{d \rightarrow \infty} \frac{Am+B}{A} \cdot \log(d)\,.
\end{equation}
\end{theorem}

According to Theorem~\ref{PACD}, all networks from the T-subclass with $A<\frac 1 2$ are assortative. However, asymptotically $\E d_{nn}(d)$ increases slowly, as $\log(d)$, unlike $d^\nu$ often observed in real-world networks. We discuss this in more details in Section~\ref{experiments}.

\subsection{Proofs}\label{Proves}

\subsubsection{Proof of Theorem \ref{PASexp}}

In the proof we use the notation $\theta(\cdot)$ for error terms.
By $\theta(X)$ we denote an arbitrary function such that $|\theta(X)| < X$.

We need the following auxiliary theorem.

\begin{theorem}\label{SNDlemm}
Let $W_n$ be the sum of the squares of the degrees of all vertices in a model from the PA-class with $A<\frac 1 2$. Then for any $\epsilon>0$
$$\E W_n = \frac{m}{1-2A} \left(m + 4B + 1 \right) n + O(n^{2A+\epsilon}).
$$
\end{theorem}

This statement is mentioned in \cite{GPA} and it can be proven by induction on $n$. We omit the proof since it is straightforward and its main idea is similar to the idea of the proof of Theorem~\ref{PASexp} described below. 



\vspace{0.3cm}

We prove Theorem \ref{PASexp} by induction on $d$ and for each $d$ we use induction on $n$. Let us prove that 
\begin{equation}\label{eq:induction}
\E S_n(d) = M(d) \, \left(n+\theta \left(C n^{2A+\varepsilon} d^{\xi}\right) \right)
\end{equation}
for some constant $C>0$.

\paragraph{Base case.}

First, for all $d$ let us obtain the basis. Namely, we prove the theorem for all $n \le Q d^2$, where $Q$ is some constant which will be specified further in the proof. Note that $N_n(d) = O\left(\frac{n}{d}\right)$, since we have $O(n)$ edges in any graph $G_m^n$. Also, by the same reason, the sum of the degrees of the neighbors of each vertex equals $O(n)$. Therefore, we have $S_n(d) = O\left(\frac{n^2}{d}\right)$ and, for $n \le Q d^2$ and $\xi \ge 3+1/A-4A$, we have $\E S_n(d) = M(d)O\left(n^{2A+\varepsilon} d^{\xi}\right)$ and so Equation~\eqref{eq:induction} holds. 

\paragraph{Inductive step for $d=m$, $n > Q \, d^2$.}

At each step of the process we add a vertex $n+1$ and $m$ edges. The following events may affect $S_n(m)$.

\begin{enumerate}
\item At least one edge hits a vertex of degree $m$, then $S_n(m)$ is decreased by the sum of the degrees of the neighbors of this vertex. This happens with probability $\frac{A \, m + B}{n} + O\left(\frac{1}{n^2} \right)$. Summing over all vertices of degree $m$ we obtain that $\E S_n(m)$ is decreased by
$
\left(\frac{Am+B}{n} + O\left(\frac{1}{n^2} \right) \right) \cdot \E S_n(m)\,.
$
\item Exactly one edge hits a neighbor of a vertex of degree $m$ and no edges hit the vertex itself, then $S_n(m)$ is increased by $1$. The probability to hit a neighbor is $\frac{Ad_i + B}{n} + O\left(\frac{d_i^2}{n^2} \right)$, where $d_i$ is the degree of this neighbor. We have to subtract the probability to hit both a vertex of degree $m$ and its neighbor which is $\frac{D}{mn} + O\left(\frac{md_i}{n^2} \right)$. Summing over all neighbors of all vertices of degree $m$, we obtain that $\E S_n(m)$ is increased by:
\begin{multline*}
\frac{A \E S_n(m)}{n} + \frac{B-D/m}{n} m \E N_n(m) + O\left(\frac{\E \sum_{\substack{i: i\text{ is a neighbor } \\ \text{ of a vertex of degree  $m$}}} d_i^2}{n^2} \right) \\ = \frac{A \E S_n(m)}{n} + \frac{B-D/m}{n} m \E N_n(m) + O\left(\frac{\max\{n, n^{3A}\}}{n^2} \right).
\end{multline*}

Here we used the fact that:
$$
\E \bigg(\sum_{\substack{i: i\text{ is a neighbor } \\ \text{ of a vertex of degree  $m$}}} d_i^2 \bigg) \le \E \left(\sum_{i \in V(G_n^m)} {d_i^3} \right)
= O\left(\max\{n, n^{3A}\}\right).
$$
\item If $i>1$ edges hit a neighbor $j$ of a vertex of degree $m$, which happens with probability $O\left(\frac{d_j^2}{n^2} \right)$, and no edges hit the vertex itself, then $S_n(m)$ is increased by $i$. Reasoning as above, we obtain that $\E S_n(m)$ is increased by~$O\left(\frac{\max\{n, n^{3A}\}}{n^2} \right).$
\item The vertex $n+1$ connects to some vertices, so $S_n(m)$ is increased by the sum of the degrees of these vertices. The probability to hit a vertex of degree $d_i$ is $\frac{Ad_i+B}{n} + O\left(\frac{d_i^2}{n^2} \right)$ and after that this vertex will have a degree $d_i+1$. Summing over $i$ we obtain that $\E S_n(m)$ is increased by:
\begin{multline}\label{eq_base:4}
\E \sum_{i \in V(G_n^m)} {(d_i+1) \left(\frac{Ad_i+B}{n} + O\left(\frac{d_i^2}{n^2} \right) \right)}   = \frac{A}{n} \E W_n + (2B+1)m + O\left(\frac{\max\{n, n^{3A}\}}{n^2} \right).
\end{multline}
\end{enumerate}

Combining all the cases considered above, we get
\begin{multline}\label{Srec}
\E S_{n+1}(m) = \E S_n(m) - \left(\frac{Am+B}{n} +O\left(\frac{1}{n^2}\right) \right) \E S_n(m) + \frac{A \E S_n(m)}{n} \\ + \frac{B-D/m}{n} m \E N_n(m) + \frac{A}{n} \E W_n + (2B+1)m + O\left(\frac{\max\{n, n^{3A}\}}{n^2} \right)\,.
\end{multline}

We are ready to show the inductive step. Assume that \eqref{eq:induction} holds for all $S_i(m)$ with $i \le n$ and let us show this equality for $i=n+1$. 
Using~\eqref{Srec} we get
\begin{multline}\label{eq:ind_m}
\E S_{n+1}(m) = \left(1 - \frac{A(m-1)+B}{n} +O\left(\frac{1}{n^2}\right) \right) M(m) \left(n+\theta \left(C n^{2A+\varepsilon}m^\xi \right) \right) \\ + \frac{B-D/m}{n}   m \cdot c(m,m) \left(n+O(1) \right) + \frac{A}{n} \cdot \frac{m}{1-2A} \left(m + 4B + 1 \right) n  \\ + O\left(n^{2A-1+\epsilon}\right) + (2B+1)m +  O\left(\frac{\max\{n, n^{3A}\}}{n^2} \right)\,.
\end{multline}
Here we use that $\E W_n = \frac{m}{1-2A} \left(m + 4B + 1 \right) n + O(n^{2A+\epsilon})$ and take $\epsilon < \varepsilon$.

According to \eqref{eq:M} for $d=m$:
\begin{equation}\label{M_m}
M(m) = \frac{m \cdot c(m,m)}{A(m-1)+B+1}  \left(B - \frac{D}{m} + \frac{(A(m-1)+2B+1)\cdot (Am+B+1)}{1-2A} \right).
\end{equation}
Combining~\eqref{eq:ind_m}, \eqref{M_m} and the fact that $c(m,m) = 1 / (Am+B+1)$ (see Equation~\eqref{Constant}) we get
\begin{multline*}
\E S_{n+1}(m) = M(m) (n+1) + \left(1 - \frac{A(m-1)+B}{n}\right) M(m) \theta \left(C n^{2A+\varepsilon}m^\xi \right) \\+ O\left(C n^{2A-2+\varepsilon} \right) + O(n^{2A-1+\epsilon})\,.
\end{multline*}

To complete the proof for $d=m$ we have to show that the obtained error term is not greater than $C M(m)  (n+1)^{2A+\varepsilon}$ for some large enough $C$:
\begin{multline*}
C M(m) (n+1)^{2A+\varepsilon}m^\xi \ge
\left(1 - \frac{A(m-1)+B}{n}\right) M(m) C n^{2A+\varepsilon}m^\xi 
 + O\left(C n^{2A-2+\varepsilon} \right) + O(n^{2A-1+\epsilon}).
\end{multline*}
This inequality holds for large enough $C$ and $n > Q m^2$ with some $Q$. This completes the inductive step for $d=m$.

\paragraph{Inductive step for $d>m$, $n > Q \, d^2$.}
Similarly to the previous case, at each step $n+1$ the following events may affect $S_n(d)$. 
\begin{enumerate}
\item At least one edge hits a vertex of degree $d$. In this case, $\E S_n(d)$ is decreased by
$\left(\frac{Ad+B}{n} + O\left(\frac{d^2}{n^2} \right) \right) \cdot \E S_n(d)$.
\item One edge hits a vertex of degree $d-1$,
so $S_n(d)$ is increased by the sum of the degrees of the neighbors of this vertex plus the degree of the new vertex. We get
\begin{equation*}
\left(\frac{A(d-1)+B}{n} + O\left(\frac{d^2}{n^2} \right) \right) \cdot \left(\E S_n(d-1) + m\cdot \E N_n(d-1)\right)\,.
\end{equation*}
Taking into account the case when, in addition, exactly one edge hits a neighbor of this vertex, we get that $\E S_n(d)$ is additionally increased by:
\begin{equation*}\label{eq:3}
(d-1)\E N_n(d-1) \cdot \frac{D}{mn} + O\left(\frac{(d-1) \E S_n(d-1)}{n^2} \right).
\end{equation*}
\item Exactly one edge hits a neighbor of a vertex of degree $d$ and no edges hit the vertex itself. In this case, $\E S_n(d)$ is increased by:
$$
\frac{A \E S_n(d)}{n} + \frac{B-D/m}{n} d\, \E N_n(d) + O\left(\frac{\max\{n, n^{3A}\}}{n^2} \right) + O\left(\frac{d \cdot \E S_n(d)}{n^2} \right).
$$
\item All the cases with multiple edges affect $\E S_n(d)$ by:
\begin{equation}\label{eq:4}
O\left(\frac{\max\{n, n^{3A}\}}{n^2} \right) + O\left(\frac{d^2}{n^2} \right) \E S_n(d) + O\left(\frac{d^3}{n^2} \right) \E N_n(d).
\end{equation}
\end{enumerate}

Combining all the cases considered above, we get
\begin{multline}\label{eq:Trec}
\E S_{n+1}(d)
= \E S_n(d) \left[1 - \frac{A(d-1)+B}{n} \right] + \frac{A(d-1)+B}{n} \cdot \E S_n(d-1) \\ + \left(\frac{D(d-1)}{mn} + m\frac{A(d-1)+B}{n}\right) \E N_n(d-1) + \frac{(B-D/m)d}{n} \E N_n(d) \\ + O\left(\frac{d^2}{n^2} \right) \E S_n(d)  + O\left(\frac{d^3}{n^2} \right)\E N_n(d) + O\left(\frac{\max\{n, n^{3A}\}}{n^2} \right)\,.
\end{multline}

Now let us show the inductive step. Assume that \eqref{eq:induction} holds for all $S_i(\tilde d)$ with $\tilde d < d$ and all $i$ and with $\tilde d = d$
and $i < n+1$.
Then
\begin{multline*}
\E S_{n+1}(d) = \left(1-\frac{A(d-1)+B}{n} \right)M(d)\left(n+ \theta \left(C n^{2A+\varepsilon} d^\xi\right) \right) \\ + \frac{A(d-1)+B}{n} M(d-1)\left(n+\theta \left(C n^{2A+\varepsilon} (d-1)^{\xi} \right) \right) \\ + \left(\frac{D(d-1)}{mn} + m\frac{A(d-1)+B}{n}\right) c(m,d-1) \left(n+O\left(d^{2+\frac{1}{A}}\right) \right) \\ + \frac{(B-D/m)d}{n} c(m,d) \left(n+O\left(d^{2+\frac{1}{A}} \right) \right) + O\left(\frac{d^2}{n^2} \right) M(d) \left(n+ \theta \left(C n^{2A+\varepsilon} d^\xi \right) \right) \\ + O\left(\frac{d^3}{n^2} \right) c(m,d) \left(n+O\left(d^{2+\frac{1}{A}}\right) \right) + O\left(\frac{\max\{n, n^{3A}\}}{n^2}\right)\,.
\end{multline*}

Note that, according to \eqref{eq:M},
\begin{multline}\label{eq:M(d)}
M(d) = \frac{A(d-1)+B}{A(d-1)+B+1} M(d-1) + \frac{(B-D/m)d}{A(d-1)+B+1} c(m,d) \\ + \frac{\left(\frac{D}{m}+Am\right)(d-1)+Bm}{A(d-1)+B+1} c(m,d-1).
\end{multline}
Therefore, we obtain:

\begin{multline*}
\E S_{n+1}(d)
= M(d)(n+1) + \left(1-\frac{A(d-1)+B}{n}\right) M(d) \, \theta \left(C n^{2A+\varepsilon} d^{\xi} \right) \\ + \frac{A(d-1)+B}{n}M(d-1)\,\theta \left(C n^{2A+\varepsilon} (d-1)^{\xi}\right) + O \left(C \frac{d^{\xi+2-\frac 1 A} \log(d) \cdot n^{2A+\varepsilon}}{n^2} \right) \\ + O\left(\frac{\max\{n, n^{3A}\}}{n^2} \right)  + O\left(\frac{d^{2}}{n} \right) + O\left(\frac{d^4}{n^2}\right).
\end{multline*}

It remains to prove that for some large enough $C$
\begin{multline}\label{eq:remains}
C M(d) \cdot  (n+1)^{2A+\varepsilon}d^{\xi}   \ge C M(d) \cdot  n^{2A+\varepsilon}d^{\xi} - C M(d) \, \left(\frac{A(d-1)+B}{n} \right) \cdot n^{2A+\varepsilon}  d^{\xi}  \\ +  C M(d-1) \,\left(\frac{A(d-1)+B}{n} \right) \,  \cdot n^{2A+\varepsilon} (d-1)^{\xi}  + O \left(C \frac{d^{\xi+2-\frac 1 A} \log(d) \cdot n^{2A+\varepsilon}}{n^2} \right) \\ + O\left(\frac{\max\{n, n^{3A}\}}{n^2} \right) + O\left(\frac{d^{2}}{n} \right)+ O\left(\frac{d^4}{n^2}\right).
\end{multline}

First, note that
\begin{multline*}
C M(d) \cdot  (n+1)^{2A+\varepsilon}d^{\xi} - C M(d) \cdot  n^{2A+\varepsilon}d^{\xi} =
  C M(d) \cdot n^{2A+\varepsilon} \cdot d^{\xi} \left(\frac{2A+\varepsilon}{n} + O\left(\frac{1}{n^2}\right) \right).
\end{multline*}

Second, one can show that
$$
C M(d) \, \left(\frac{A(d-1)+B}{n} \right) \, d^{\xi} - C M(d-1) \, \left(\frac{A(d-1)+B}{n} \right)\,  (d-1)^{\xi} \ge 0
$$
using Equation~\eqref{eq:M(d)}, the inequality $(1-\frac{1}{d})^{-\xi} \ge 1 + \frac{\xi}{d}$, and the fact that $\xi \ge \frac{2}{1-A}$.


Therefore, Equation \eqref{eq:remains} becomes:

\begin{multline*}
  C M(d) \cdot n^{2A+\varepsilon} \cdot d^{\xi} \left(\frac{2A+\varepsilon}{n} + O\left(\frac{1}{n^2}\right) \right) \ge
O\left(C \frac{d^{\xi+2-\frac{1}{A}} \log(d) \cdot n^{2A+\varepsilon}}{n^2} \right) \\ + O\left(\frac{\max\{n, n^{3A}\}}{n^2} \right) + O\left(\frac{d^{2}}{n} \right) + O\left(\frac{d^4}{n^2}\right).
\end{multline*}



It is easy to see that for some large enough $C$ and for $n \ge Q \cdot d^2$ (for some large enough $Q$) this inequality is satisfied. This concludes the inductive step for $d>m$ and also the proof of the main part of the theorem (Equation~\eqref{eq:main}).

Now, let us show why approximation~\eqref{eq:M_limit} holds.
First, we estimate $\sum_{i=m+1}^{d} {Y(i)}$ in Equation~\eqref{eq:M}:
\begin{multline*}
\sum_{i=m+1}^{d} {Y(i)} 
= \sum_{i=m+1}^{d} \frac{1}{A(i-1)+B+1} \left(\frac{(B-D/m)i}{Ai+B+1}+\frac{(D/m)\cdot(i-1)}{A(i-1)+B} +m \right) \\ \widesim{d \rightarrow \infty} \frac{Am+B}{A^2} \log(d).
\end{multline*}
Note that $X$ is a constant.
Also, recall that 
$c(m,d) \widesim{d \rightarrow \infty} \frac{\G\left(m + \frac{B+1}{A}\right)d^{-1-\frac{1}{A}}}{A\, \G\left(m + \frac{B}{A}\right)}$ (see Eqution~\eqref{Constant}).
Finally, 
\begin{multline*}
M(d) = \left(Ad+B+1\right) \left(\frac{X}{Am+B+1} + \sum_{i=m+1}^{d} {Y(i)}\right) \cdot c(m,d) \\
 \widesim{d \rightarrow \infty} Ad \cdot \frac{Am+B}{A^2} \log(d) \cdot \frac{\G\left(m + \frac{B+1}{A}\right)d^{-1-\frac{1}{A}}}{A\, \G\left(m + \frac{B}{A}\right)} =
\frac{Am+B}{A^2} \cdot \frac{\G\left(m + \frac{B+1}{A}\right)}{\G\left(m + \frac{B}{A}\right)} \cdot \log(d) \cdot d^{-\frac{1}{A}}.
\end{multline*}
This completes the proof.

\subsubsection{Proof of Theorem \ref{PACD}}

Denote by $Q$ the event $\left\{|N_n(d) - \E N_n(d)| < d \sqrt{n} \log(n) \right\}$.
According to Theorem~\ref{Concentration}, $\Prob(Q) = 1 - O\left(n^{-\log(n)}\right)$.
Let us estimate $\E d_{nn}(d)$:
$$\E d_{nn}(d) =
\E \left(\frac{S_n(d)}{d \, N_n(d)}\right) = \E \left(\frac{S_n(d)}{d \, N_n(d)}\bigg|Q\right)\Prob(Q) + \E \left(\frac{S_n(d)}{d N_n(d)}\bigg|\bar Q\right)\Prob({\bar Q}).$$

Let us estimate the first term:
\begin{multline*}
\E \left(\frac{S_n(d)}{d \, N_n(d)}\bigg|Q\right)\Prob(Q) =
\frac{\E \left(S_n(d)\big|Q\right)\Prob(Q)}{d \left( \E N_n(d) + O\left(d \sqrt{n} \log(n)\right)\right)}\\
= \frac{\E S_n(d) - \E (S_n(d) \big|\bar Q) \Prob(\bar Q)}{d \left( \E N_n(d) + O\left(d \sqrt{n} \log(n)\right)\right)}
= \frac{\E S_n(d) + O\left(n^{1-\log(n)}\right)}{d \left( \E N_n(d) + O\left(d \sqrt{n} \log(n)\right)\right)}.
\end{multline*}
Here we used that $S_n(d) = O(n)$.
The second term can be estimated as:
$$
\E \left(\frac{S_n(d)}{d\,N_n(d)}\bigg|\bar Q\right)\Prob({\bar Q})
= O\left(\frac{n}{d} \right)\Prob({\bar Q})
= O\left(\frac{n^{1-\log(n)}}{d}\right).
$$
Finally,
\begin{multline*}
\E d_{nn}(d) = \frac{M(d)\left(n + O \left( n^{2A+\varepsilon} \cdot d^{\xi}\right)\right) + O\left(n^{1-\log(n)}\right)}{d \left( c(m,d)\left(n + O \left( d^{2+\frac{1}{A}}\right)\right) + O\left(d \sqrt{n} \log(n)\right)\right)} + O\left(\frac{n^{1-\log(n)}}{d}\right) \\ =
\frac{M(d)}{d \, c(m,d)}\left(1 + O\left(\frac{n^{2A+\varepsilon} \cdot d^{\xi}}{n} + \frac{d^{2+\frac 1 A} \log{(n)}}{\sqrt{n}} \right) \right).
\end{multline*}

\subsection{Hypotheses ($A \ge 1/2$)}\label{hypothesis}

Note that the restriction $A<\frac{1}{2}$ is essential and for $A \ge \frac{1}{2}$ the result is expected to be completely different.
For example, similarly to the configuration model discussed in~\cite{ANND}, here we expect $d_{nn}(d)$ to scale with $n$.

In the proof of Theorem~\ref{PASexp} we first analyze $\E S_n(m)$ and we have to estimate the expected sum of the degrees of the neighbors of a new vertex $n+1$ (see Equation~\eqref{eq_base:4}).  If $A \ge \frac{1}{2}$, then the term $\frac{A}{n}\E W_n$ in Equation~\eqref{eq_base:4} grows with $n$ (the behavior of $\E W_n$ for $A\ge \frac 1 2$ is discussed, e.g., in~\cite{GPA}). This leads to the fact that $\E S_n(m)$ grows with $n$ faster than linearly.
Also, it can be shown that 
the initial configuration $G_m^{n_0}$ for any constant $n_0$ affects $\E W_n$ by at least constant multiplicative factor (this can be shown by following the proof of Theorem~\ref{PASexp}), $\E W_n$, in tern, affects $\E S_n(m)$, and so all $S_n(d)$ are affected.
This means that if $A \ge \frac 1 2$, then precise asymptotics for $\E d_{nn}(d)$, similar to the one discussed in Theorem~ \ref{PACD}, cannot be obtained for the whole T-subclass, as this asymptotics depends on the initial configuration $G_m^{n_0}$.
Also, even if the initial configuration $G_m^{n_0}$ is fixed, we expect that (in the case of infinite variance of the degree distribution) $d_{nn}(d)$
will probably not converge as $n \to \infty$. For a deeper discussion of this phenomenon see~\cite{ANND}.
Therefore, in this section we only propose some hypotheses which we then test by simulations.


Consider the case $A > \frac 1 2$. Similarly to Theorem~\ref{SNDlemm}, one can see that $W_n$ asymptotically behaves as $n^{2A}$ (see also~\cite{GPA}). Therefore, let us assume that $\E W_n = C_1 \cdot n^{2A}$ for some constant $C_1$. We can now approximate $\E S_n(d)$ by  $M(d) \cdot n^{2A}$ with some function $M(d)$ which we will compute now. From \eqref{Srec} we get (omitting all error terms):
$$
M(m) (n+1)^{2A} = M(m) n^{2A} - \frac{A(m-1)+B}{n} M(m) n^{2A} + \frac{A}{n} C_1 n^{2A},
$$
so we can approximate $M(m)$ by
$$
M(m) = \frac{AC_1}{A(m+1)+B}.
$$
Similarly, for $d>m$ we have
$$
M(d) (n+1)^{2A} = M(d) n^{2A} - \frac{A(d-1)+B}{n} n^{2A} \left(M(d-1) - M(d) \right),
$$
and we can get the following approximation
$$
M(d) = M(d-1) \cdot \frac{A(d-1)+B}{A(d+1)+B} = M(m) \prod_{i=1}^{d-m} {\frac{A(d-i)+B}{A(d+2-i)+B}} \\ = \frac{AC_1(Am+B)}{(Ad+B)(A(d+1)+B)}.
$$
Finally, in the case $A > \frac{1}{2}$ we have 
\begin{align}\label{hyp1}
\E d_{nn}(d) \sim \frac{AC_1(Am+B)}{(Ad+B)(A(d+1)+B) \,d \, c(m,d)} \cdot n^{2A-1} \\ \sim C_1(Am+B) \frac{\G\left(m+\frac{B}{A}\right)}{\G\left(m+\frac{B+1}{A}\right)} \cdot d^{\frac{1}{A}-2} n^{2A-1}. \nonumber
\end{align}

By similar reasoning we can obtain an approximation of $\E d_{nn}(d)$ for $A = \frac 1 2$. In this case, we assume that $\E W_n = C_2 \, n \log(n)$ and approximate $\E S_n(d)$ by  $M(d) \cdot \log(n)$ with some $M(d)$. Finally, substituting everywhere $A = \frac 1 2$ and $B = 0$, we get
\begin{align}\label{hyp2}
\E d_{nn}(d) \sim \frac{C_2\,(d+2)}{2\,d\,(m+1)} \cdot \log(n) \\ \sim \frac {C_2} {2\,(m+1)}\cdot \log(n). \nonumber
\end{align}

In the next section, we conduct several experiments and illustrate, in particular, the above hypotheses. 

\section{Experiments}\label{experiments}

In this section, we first illustrate our theoretical results from Section~\ref{main_results} and then our hypotheses discussed in Section~\ref{hypothesis}. For these purposes, we generate several graphs using a three-parameter model from the family of polynomial graph models defined in~\cite{GPA}. This model belongs to the T-subclass and by varying the parameters we can analyze the effect of $A$ (or, equivalently, $\gamma$) and $D$ on $d_{nn}(d)$. Detailed graph generation process is described in~\cite{GPA}. 

\subsection{$A<1/2$}

\begin{figure*}[t]
        \centering
        \begin{subfigure}[b]{0.48\textwidth}
            \centering
            \includegraphics[width=\textwidth]{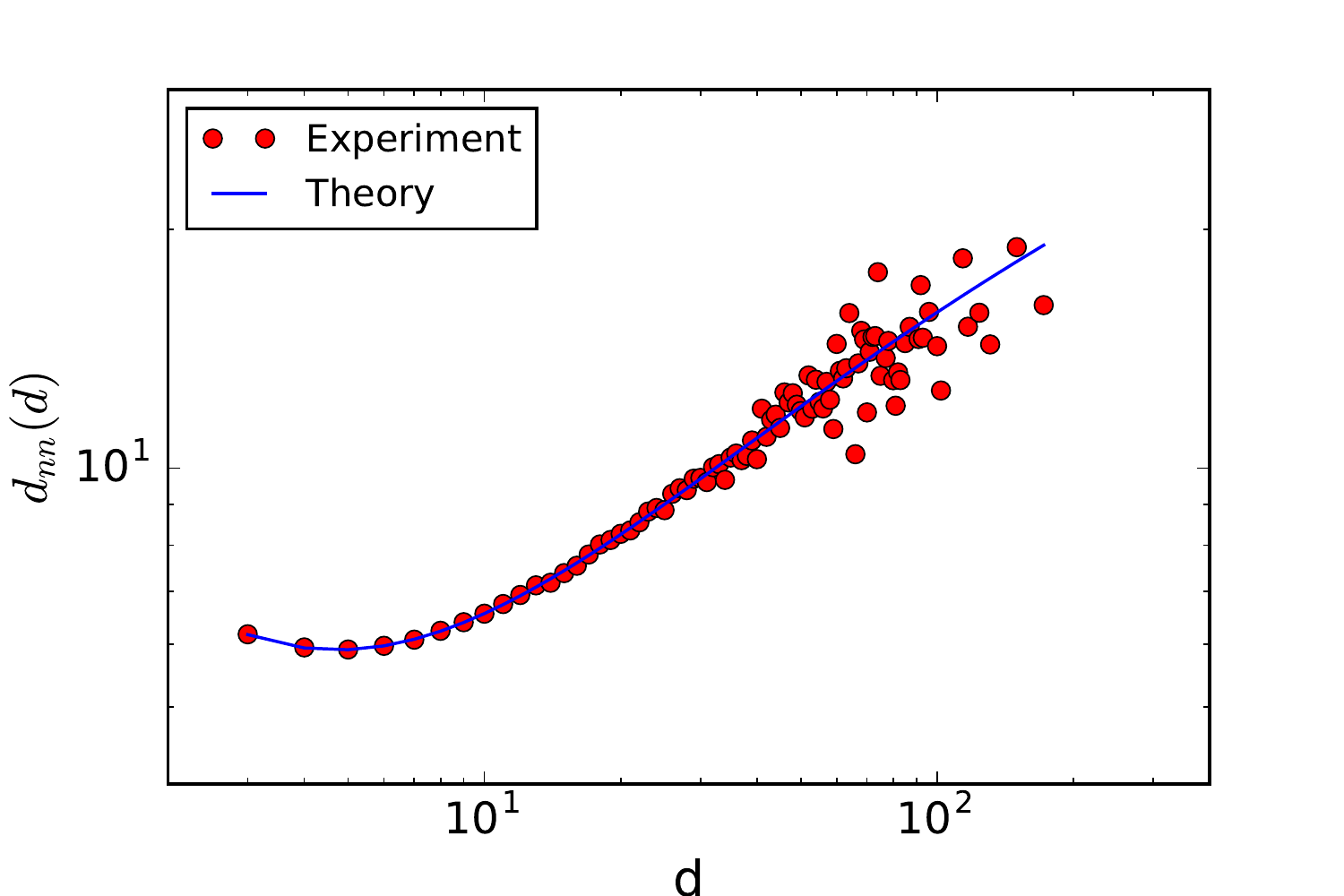}
            \caption{$A=0.2$}
        \end{subfigure}
        \begin{subfigure}[b]{0.48\textwidth}
            \centering
            \includegraphics[width=\textwidth]{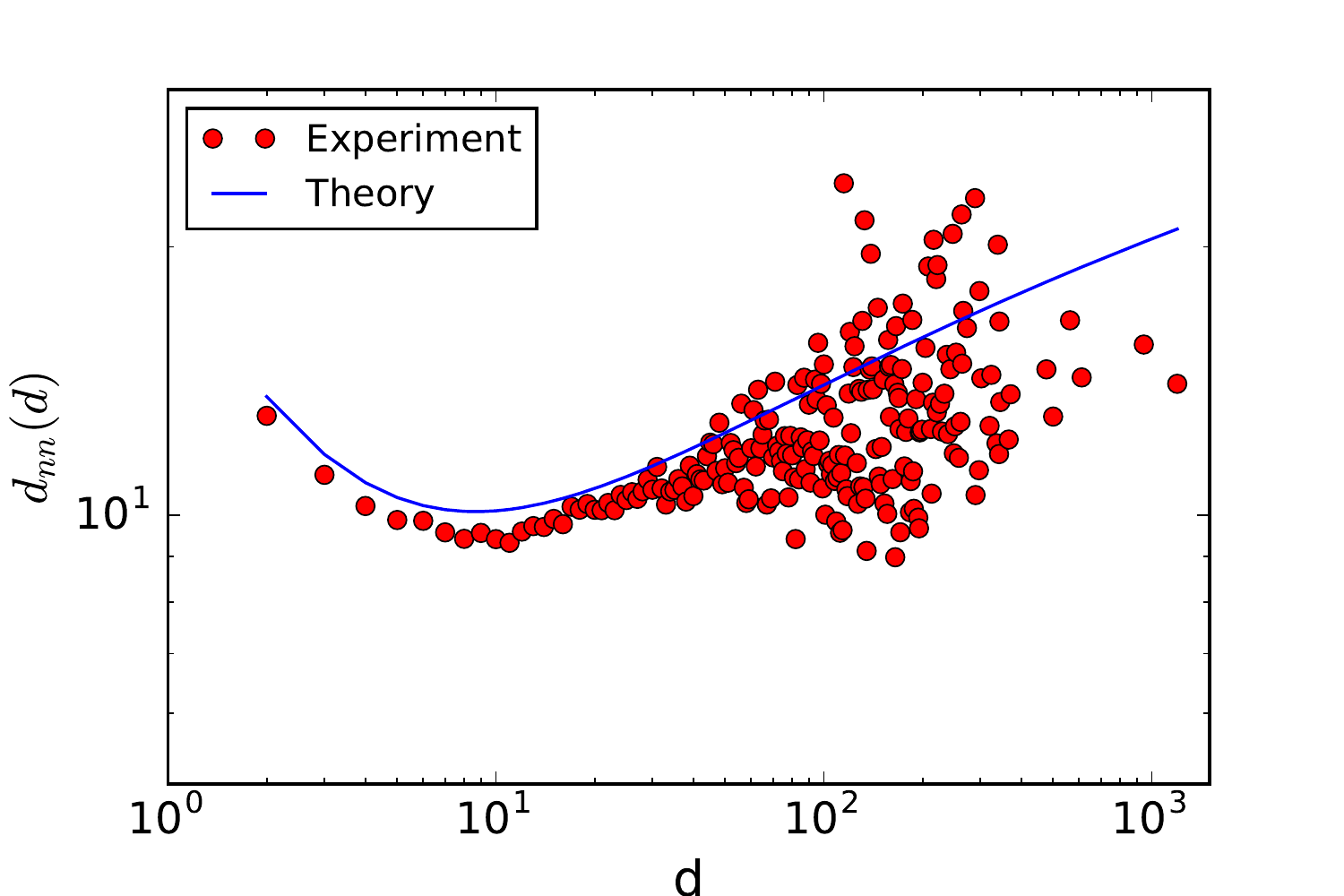}
            \caption{$A=0.4$}
        \end{subfigure}
        \caption{The behavior of $d_{nn}(d)$ for $A< 1/2$}
        \label{fig:C_D_graphs}
\end{figure*}

\begin{figure}[t]
	\begin{center}
		\includegraphics[height = 5cm]{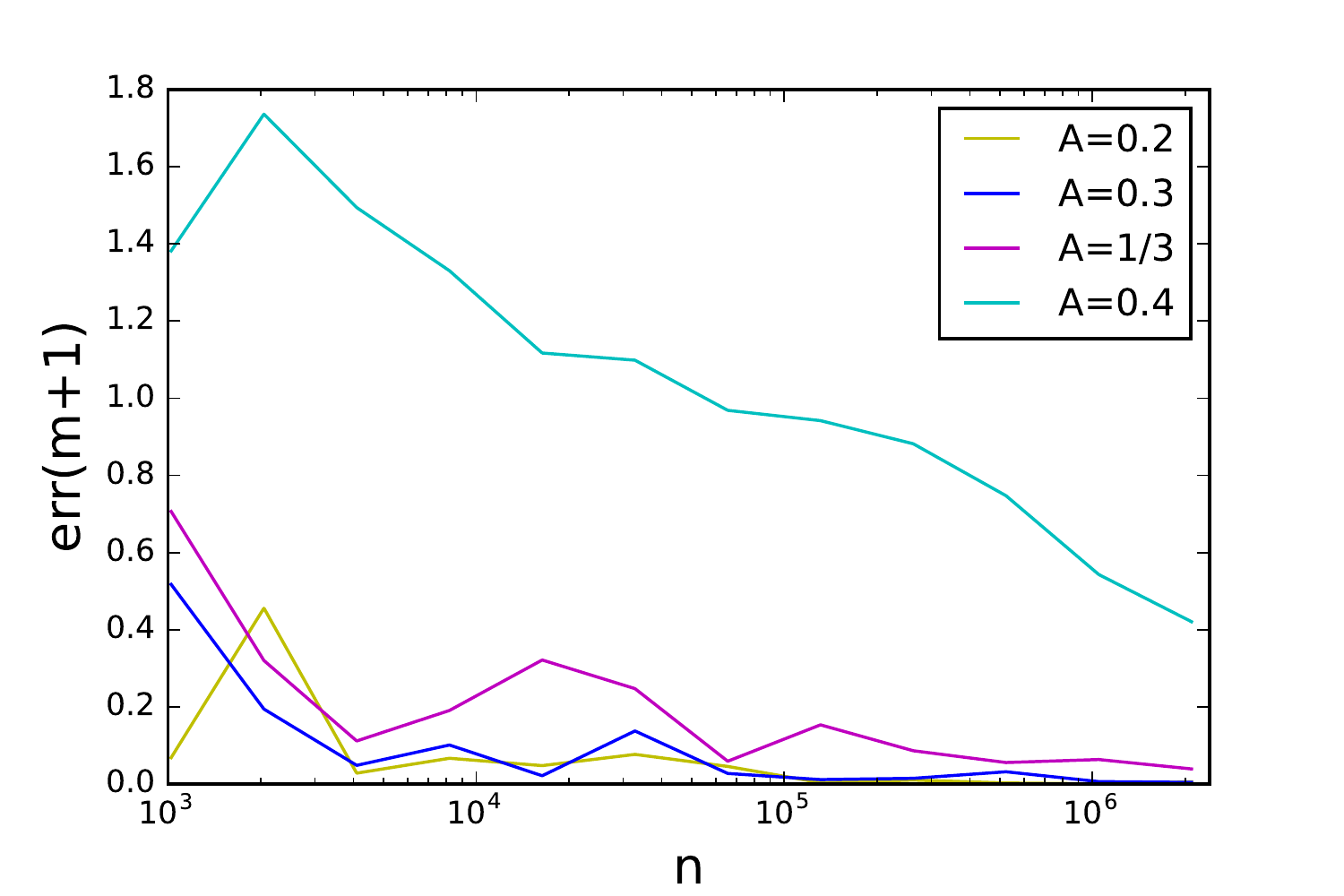}
	\end{center}
	\caption{Convergence of $d_{nn}(m+1)$ to its approximation}
	\label{fig:err}
\end{figure}

First, let us illustrate our main result for $\E d_{nn}(d)$ (see Theorem~\ref{PACD}).
We generated several polynomial graphs and compared the obtained values of $d_{nn}(d)$ with their theoretical approximation $\frac{M(d)}{d\cdot c(m,d)}$. We took $n=10^6, m=2, D=0.3$ and considered different values of $A$. In other words, we fixed the probability of a triangle formation and varied the parameter of the power-law degree distribution. 
We noticed that for $A < \frac{1}{3}$ the theoretical value of $\E d_{nn}(d)$ is extremely close to the experiment.
However, if $\frac 1 3 < A < \frac 1 2$, then $d_{nn}(d)$ turn out to be consistently smaller than their theoretical approximation.
Figure~\ref{fig:C_D_graphs} illustrates this observation and shows $d_{nn}(d)$ for $A = 0.2$ and $A = 0.4$. We decided to deeper analyze this phenomenon and plotted the difference between theoretical and empirical values of $d_{nn}(d_0)$ for some $d_0$ and different values of $A$. We took $A\in\{0.2, 0.3, 1/3, 0.4\}$ and $d_0 = m+1 = 3$\footnote{In this and further experiments we choose $d_0 = m+1$. We want to cover at least one induction step, therefore we can take any $d_0 > m$.}. Figure~\ref{fig:err} presents the behavior of $\mathrm{err}(m+1) = \left|d_{nn}(m+1) - \frac{M(m+1)}{(m+1)\cdot c(m,m+1)}\right|$ averaged over 10 samples of graphs generated by the polynomial model with $D=0.2$ and different values of $n$. From Figure~\ref{fig:err} we see that the difference between $d_{nn}(m+1)$ and its theoretical approximation decreases as $n$ grows for all values of $A$. However, for  $A=0.4$ the converges is much slower. The possible reason for this is the error term $O\left(\frac{n^{3A}}{n^2}\right)$ appearing in the proof in the case $A > \frac 1 3$. 

\begin{figure}[t]
	\begin{center}
		\includegraphics[height = 5cm]{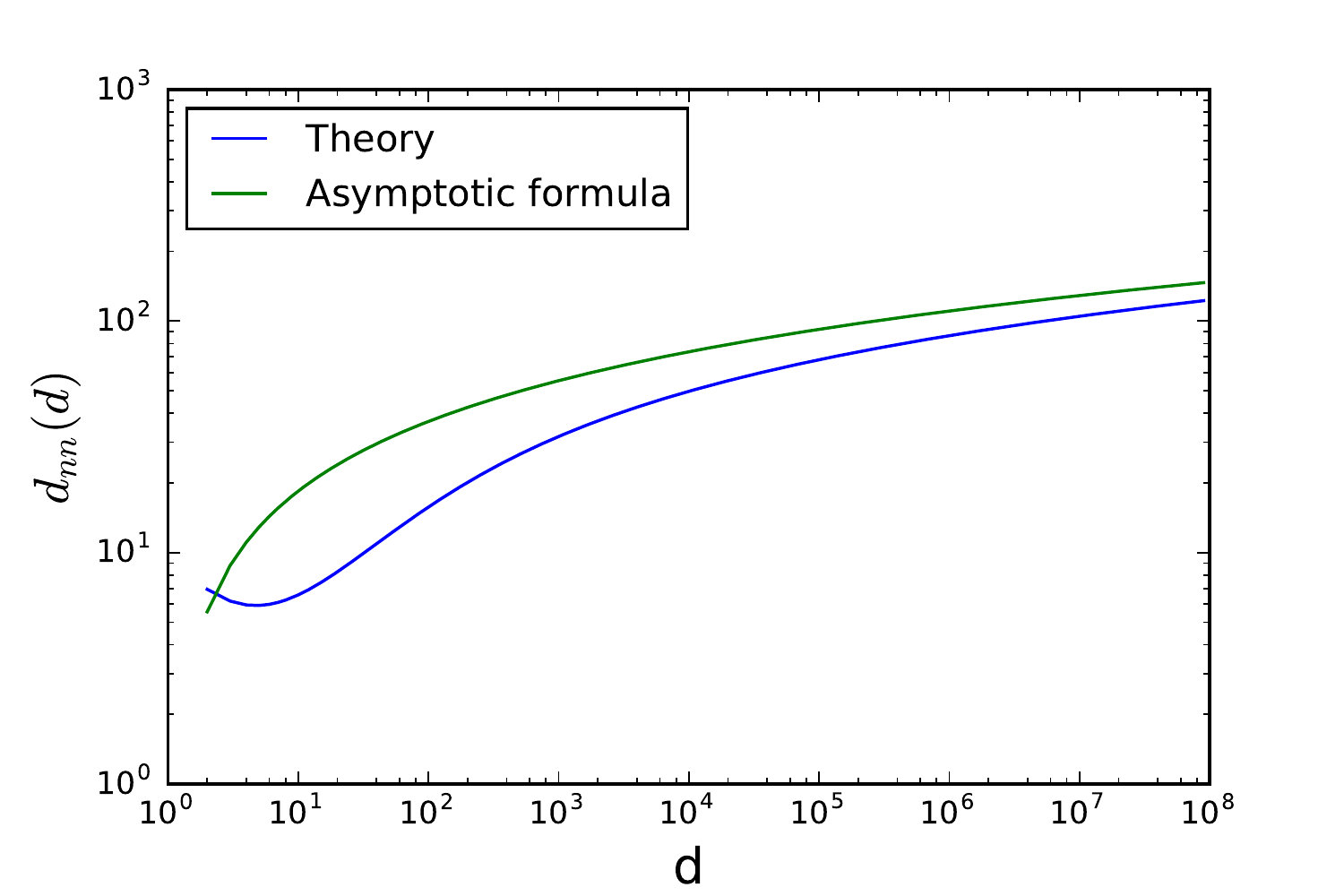}
	\end{center}
	\caption{Theoretical value of $d_{nn}(d)$ versus its asymptotic approximation}
	\label{fig:Clust_C_2_D}
\end{figure}

We also compared the theoretical value of $\E d_{nn}(d)$ (for $A=0.2$, $D=0.3$) with the asymptotic formula $\frac{Am+B}{A} \cdot \log(d)$ (see Figure~\ref{fig:Clust_C_2_D}).
Interestingly, from Figure~\ref{fig:C_D_graphs} it may seem that $d_{nn}(d)$ grows as $d^\nu$ for some $\nu$ (as it was observed in many real-world networks). However, as $d$ becomes large ($d> 10^4$), one can indeed observe the logarithmic growth.

\begin{figure}[t]
	\begin{center}
		\includegraphics[height = 5cm]{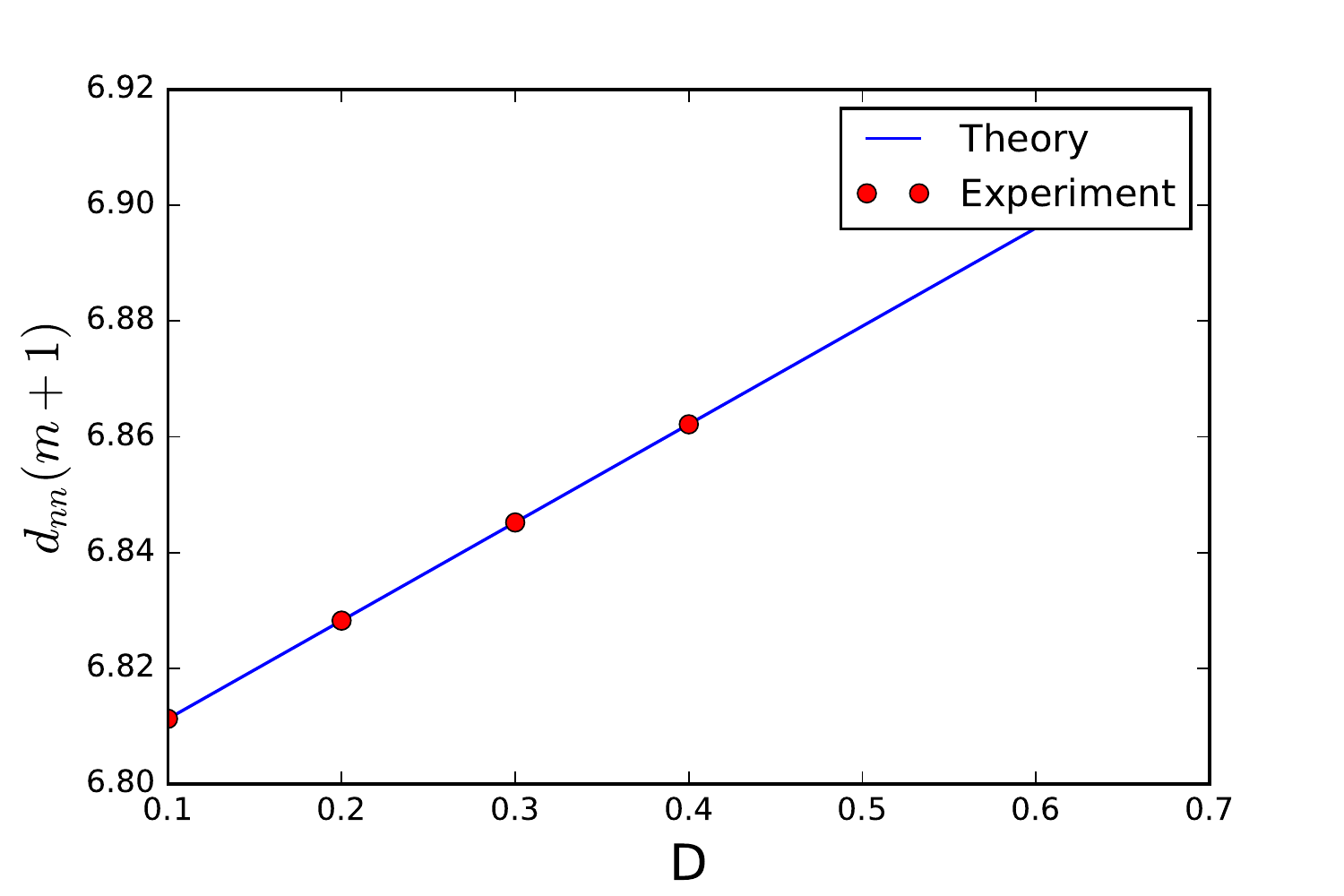}
	\end{center}
	\caption{$d_{nn}(m+1)$ versus its approximation depending on $D$}
	\label{fig:Clust_C_2_D_05_D}
\end{figure}

Finally, we looked at the behavior of $d_{nn}(d)$ and its approximation $\frac{M(d)}{d\cdot c(m,d)}$ depending on $D$. For this purpose, we plotted $d_{nn}(d_0)$ averaged over 10 samples for different values of $D$. We took $m=2$, $d_0=m+1$, $n=10^6$, and $A=0.25$ (see Figure~\ref{fig:Clust_C_2_D_05_D}). First, note that the dependence on $D$ is almost negligible. Second, $d_{nn}(m+1)$ grows linearly with $D$. Indeed, by the definition (see Theorem~\ref{PASexp}) $M(d)$ depends linearly on $D$ and is asymptotically independent of $D$ (for $n \to \infty$).

\subsection{$A \ge 1/2$}

In this section, we illustrate our hypotheses discussed in Section~\ref{hypothesis}. For this, we generated polynomial graphs with $A\in\{0.5,0.6\}$, $D=0.2$, $m=2$ and compared
 $d_{nn}(d)$ with our hypotheses~\eqref{hyp1} and~\eqref{hyp2}. 
 
\begin{figure*}[t]
        \centering
        \begin{subfigure}[b]{0.48\textwidth}
            \centering
            \includegraphics[width=\textwidth]{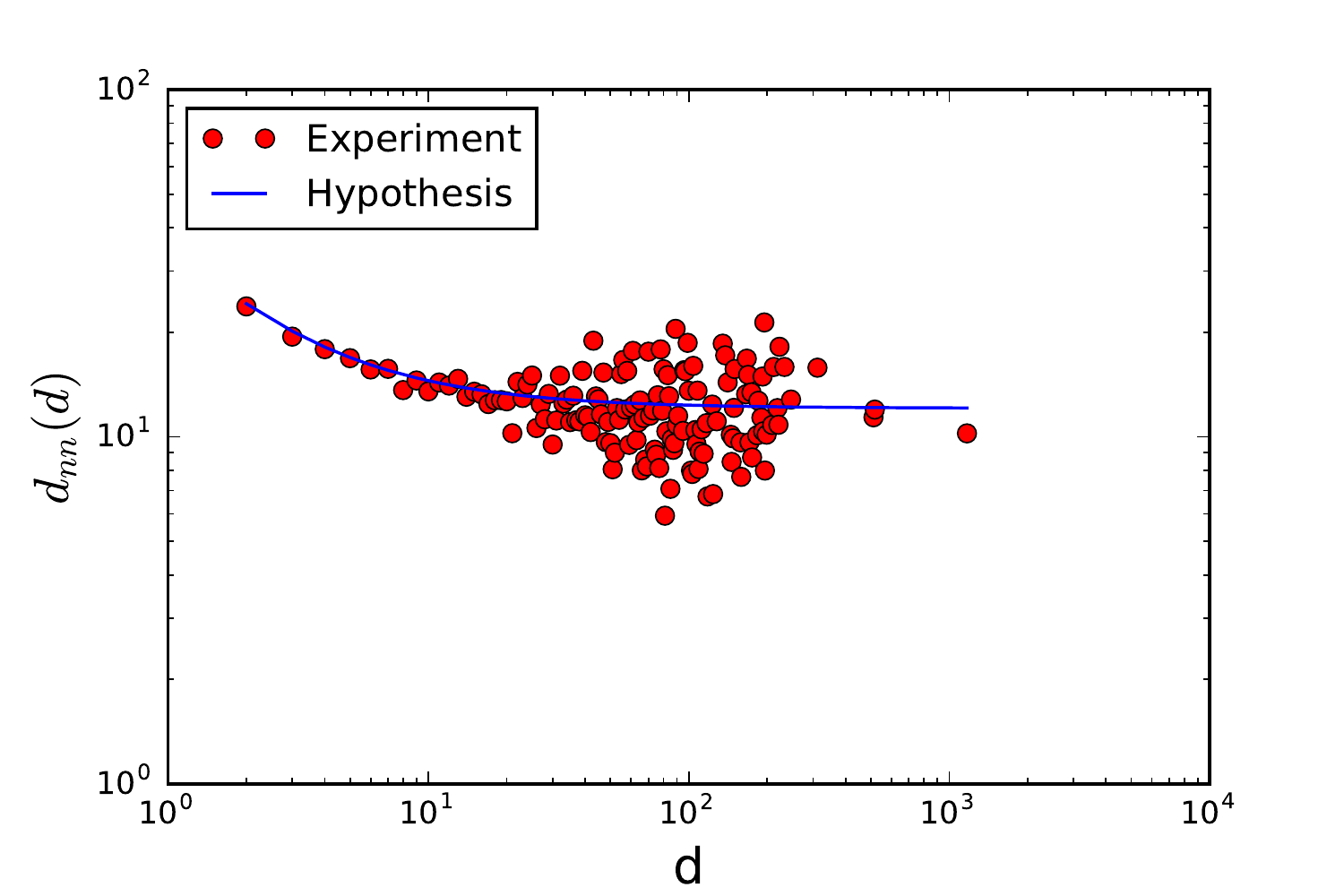}
            \caption{$A=0.5$}
            \label{fig}
        \end{subfigure}
        \hfill
        \begin{subfigure}[b]{0.48\textwidth}
            \centering
            \includegraphics[width=\textwidth]{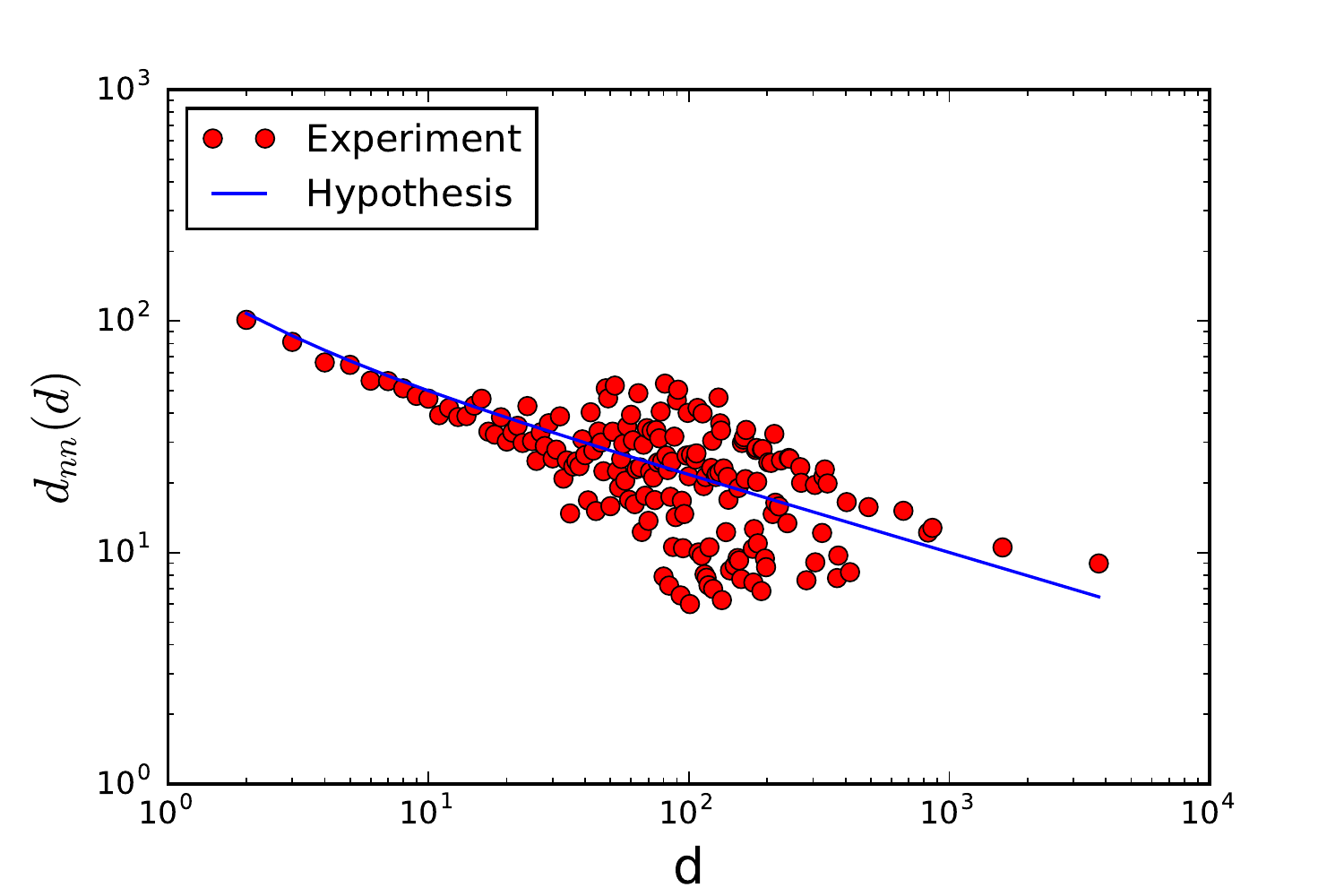}
            \caption{$A=0.6$}
	 \end{subfigure}
        \caption{$d_{nn}(d)$ versus our hypothesis for $A\ge 1/2$}\label{fig:hyp_real}
        \label{fig:hyp_d}                
\end{figure*}

Figure~\ref{fig:hyp_d} illustrates empirical values of $d_{nn}(d)$ and our approximations~\eqref{hyp1} and~\eqref{hyp2}  for $n = 10^6$.
Recall that for $A>\frac 1 2$ we expect (asymptotically) $\E d_{nn}(d) \propto n^{2A-1} d^{\frac 1 A - 2}$ and for $A = \frac 1 2$ we have $\E d_{nn}(d) \propto \log(n)$, which does not depend on $d$.
We see that for both values of $A$ our hypotheses are very close to the observed behavior.

\begin{figure*}[t]
	\centering
	\begin{subfigure}[b]{0.48\textwidth}
		\centering
		\includegraphics[width=\textwidth]{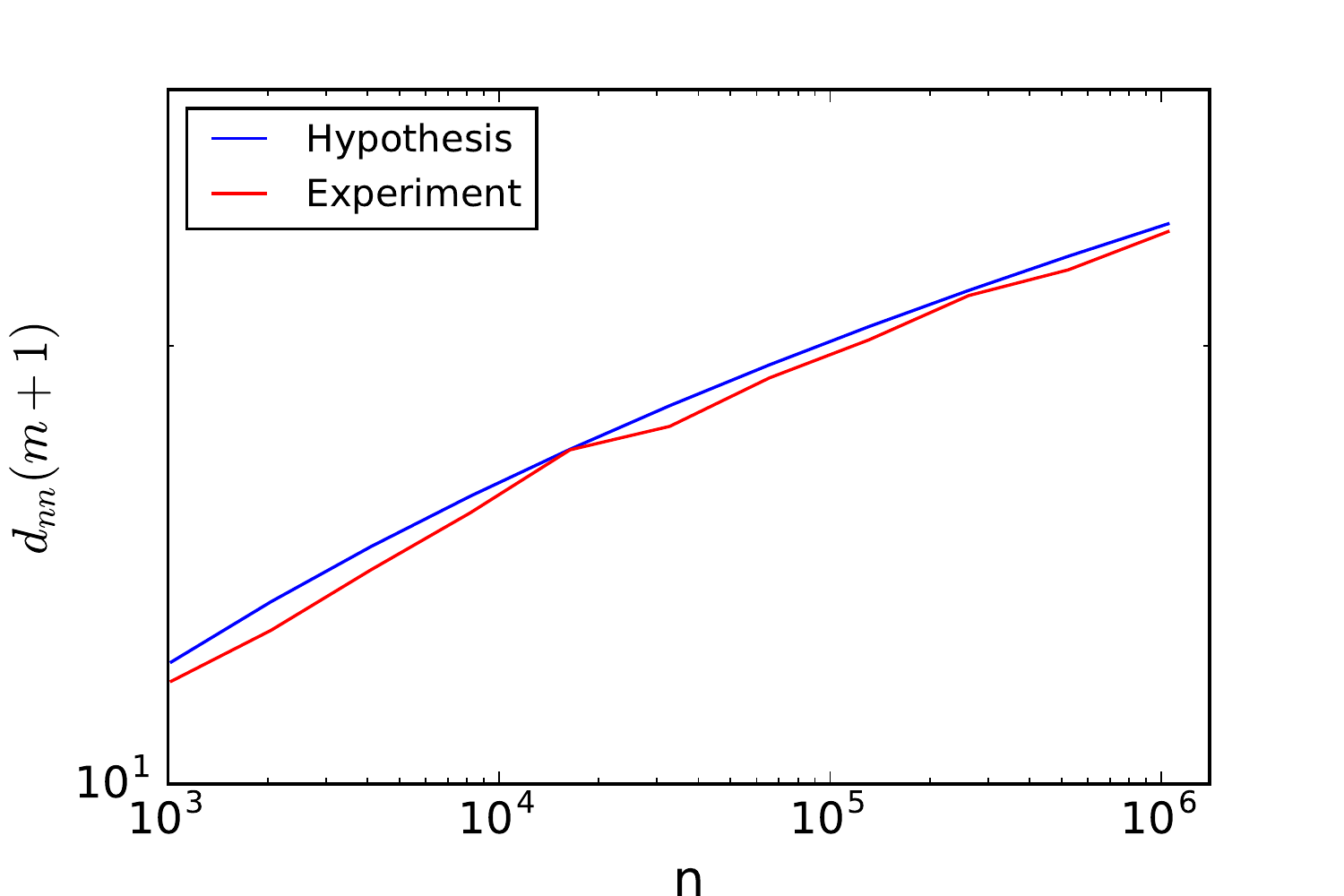}
		\caption{$A=0.5$}
	\end{subfigure}
	\hfill
	\begin{subfigure}[b]{0.48\textwidth}
		\centering
		\includegraphics[width=\textwidth]{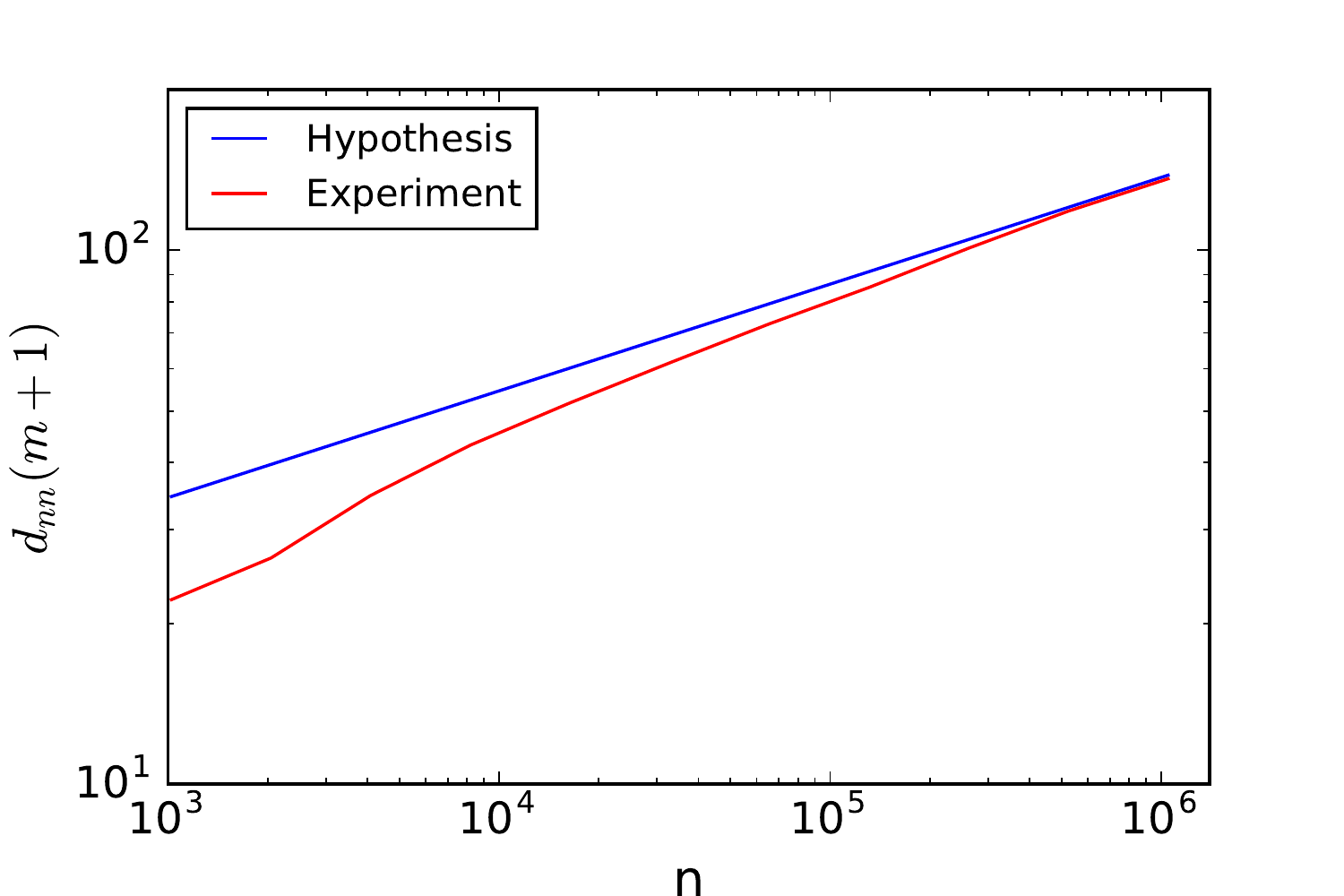}
		\caption{$A=0.6$}
	\end{subfigure}
	\caption{$d_{nn}(m+1)$ as a function of $n$ versus our hypothesis for $A\ge 1/2$}
	\label{fig:hyp_n}
\end{figure*}

As discussed above, for $A\ge \frac 1 2$ we expect $\E d_{nn}(d)$ to scale with $n$.
Therefore, we additionally analyzed the dependence of $d_{nn}(d_0)$ on $n$ for $d_0 = m+1$ (see Figure~\ref{fig:hyp_n}) and in this case we also obtained a good approximation. Note that the constants $C_1$ and $C_2$ in Equations~\eqref{hyp1} and~\eqref{hyp2} were chosen manually, but they are the same for Figures~\ref{fig:hyp_d} and~\ref{fig:hyp_n}.

Also, it worth noting that the case $A > \frac 1 2$ corresponds to many real-world networks. In~ \cite{Empirical} it was shown that the behavior of $d_{nn}(d)$ in the Buckley-Osthus model (with carefully chosen parameter which corresponds to the degree distribution with infinite variance) is very similar to the one observed for the web host graph.




\section{Conclusion and future work}\label{sec:conclusion}

In this paper, we studied the degree-degree correlations in the PA-class of models. Namely, we analyzed the behavior of the average neighbor degree $d_{nn}(d)$. For the whole PA-class of models we estimated $\E d_{nn}(d)$ for the case $A< \frac 1 2$. In particular, we proved that $\E d_{nn}(d) \widesim{d \rightarrow \infty} \frac{Am+B}{A} \cdot \log(d)$. We also discussed the case $A\ge \frac 1 2$, argued why in this case $\E d_{nn}(d)$ scales with $n$, and proposed hypotheses on its asymptotic behavior. 

There are several important directions for future research in this area. First, note that in Theorem~\ref{PACD} we analyze only the average value of $d_{nn}(d)$ and the second step is proving concentration. Usually in such problems  the Azuma-Hoeffding inequality is used (see, e.g.,~\cite{LCC}). However, in the case of $d_{nn}(d)$ we expect additional difficulties since for each $d$ a new vertex added at some step can hugely affect $d_{nn}(d)$.
Another interesting direction is proving central limit theorem for the case $A\ge \frac 1 2$, similarly to the one in~\cite{ANND}. 
More importantly, it would be interesting to analyze the average nearest neighbor rank proposed in~\cite{ANND} instead of $d_{nn}(d)$. However, we expect this to be much more difficult for the PA-class of models.


\end{document}